\documentclass[12pt]{amsart}

\usepackage{amsfonts}
\usepackage{amsmath}
\usepackage{amssymb}
\usepackage{amstext}
\usepackage{amsthm}
\usepackage{mathtools}
\usepackage{amsopn}
\usepackage{mathrsfs}
\usepackage{geometry}
\usepackage{url}
\usepackage{amsrefs}
\usepackage{hyperref}

\usepackage{booktabs} 
\usepackage{amssymb}
\usepackage{amsthm}
\usepackage{amsmath}
\usepackage{mathtools}
\usepackage{amsopn}
\usepackage{mathrsfs}
\usepackage{sansmathaccent}
\usepackage{tikz}
\usepackage{xcolor}
\usepackage{multirow}
\usepackage{booktabs}
\usepackage{xparse}
\usepackage{calc}
\usepackage{MnSymbol,wasysym}
\usepackage{pgfplots}
\usepackage{verbatim}
\usepackage{graphicx}

\usetikzlibrary{arrows}
\usepackage{pgfplots}
\usetikzlibrary{calc}
\usetikzlibrary{datavisualization.formats.functions}

\pgfplotsset{my style/.append style={axis x line=middle, axis y line=
		middle, xlabel={$x$}, ylabel={$y$}, axis equal }}

\usepackage[english]{babel}
\usepackage [autostyle, english = american]{csquotes}
\MakeOuterQuote{"}

\font\myfont=cmr12 at 6pt

\hypersetup{
	colorlinks=true,
	linkcolor=red,
	filecolor=magenta,      
	urlcolor=cyan,
}

\theoremstyle{definition}

\newtheorem{theorem}{Theorem}[section]

\newtheorem{remark}[theorem]{Remark}

\newtheorem{counter example}[theorem]{Counter Example}

\newtheorem*{choice}{Axiom of Choice}
\newtheorem*{c-choice}{Countable Axiom of Choice - CAC} 
\newtheorem*{d-choice}{Axiom of Dependent Choice - ADC}
\newtheorem*{constructive-conclusion}{Constructive-Conclusion}

\newcommand{\cc}[1]{\ignorespaces}
\allowdisplaybreaks[1]

\newcommand{\I}{\mathbf{i}}

\DeclareMathOperator\Ln{Ln}

\title[A Paradox on the Law of Excluded Middle]{A Paradox on the Law of Excluded Middle in the framework of category of set}

\author[Babak Jabbar Nezhad]{Babak Jabbar Nezhad \\ {\myfont "Dedicated to my beloved mother Simzar Hosseinzadeh who has been a source of love, wisdom and care to me, to my father Firouz who was my first math teacher, to my brother Masood for his support, emotionally, mentally and scientifically, and to my sister Neda for her love and support."}}
\address{Da\c{s} Maku, West Azerbaijan, Iran}
\email{babak.jab@gmail.com}


\date{}
\thanks{Babak Jabbar Nezhad has also published under the name Babak Jabarnejad~\cite{jabarnejad2016rees}}

\dedicatory{}

\begin{document}
\maketitle
\vspace{-7mm}
	
\begin{abstract}
In this paper, we present a paradox arising from the acceptance of the Law of Excluded Middle (LEM) within classical mathematics. Specifically, we construct a nonzero analytic function on a connected open subset of the complex plane whose zeros are not isolated. This contradicts a fundamental theorem in complex analysis, thereby revealing an inconsistency tied to LEM. Unlike traditional critiques that reject LEM in favor of intuitionistic or constructive mathematics, we argue that LEM is instrumental in discovering \textbf{relations} between objects and facts rather than the objects themselves. Since we are not always in direct attachment with objects, this relational perspective introduces \textbf{inherent uncertainty} in mathematical reasoning. Consequently, we propose that the logical framework of the world is undecidable, making contradictions possible in more complex contexts. Our findings suggest that LEM, while powerful, may not be universally reliable in all mathematical frameworks. This work has implications for foundational mathematics, particularly in relation to the limits of classical logic and the necessity of alternative logical paradigms.
\end{abstract}

\section{Introduction}
Mathematical knowledge is built upon rigorous logical structures developed over centuries. Among the most debated principles is the Law of Excluded Middle (LEM), which asserts that for any proposition $P$, either $P$ or $\neg P$ must hold. While LEM is central to classical mathematics, it is rejected in constructive mathematics, where truth is established through explicit constructions rather than logical negation~\cite{bi, td}.

The debate surrounding LEM extends beyond mathematics into philosophy, particularly in discussions by Brouwer~\cite{br}, and Dummett~\cite{du}. However, rather than taking an intuitionistic stance, we introduce a different perspective: \textbf{LEM is not merely a logical tool, but a mechanism for discovering relations between mathematical objects, rather than the objects themselves}. Since we often operate without direct access to mathematical objects, we rely on logical structures that encode relationships, introducing an unavoidable \textbf{uncertainty} into our reasoning. This suggests that the logical structure of reality itself may be undecidable, an idea with significant implications for mathematical foundations.

In this paper, we construct a concrete mathematical example that exposes a paradox inherent in LEM. Using classical complex analysis, we define a nonzero analytic function on a connected open domain whose zeros are not isolated. Since classical complex analysis dictates that such zeros must be isolated, this leads to a contradiction, challenging the robustness of LEM in mathematical reasoning. Unlike previous counterexamples that either arise within constructivism or rely on philosophical argumentation, our paradox emerges within the \textbf{category of sets}, where classical logic is considered most stable.

The paper is structured as follows: Section~\ref{back} provides background on LEM, the Axiom of Choice, and relevant results in complex analysis. Section~\ref{ex} presents our example and demonstrates the paradox. Finally, in Section~\ref{con}, we discuss the broader implications of our findings, particularly regarding the undecidability of logical frameworks.      
  
\section{A little of background}\label{back}
First of all, one should be familiar with following axioms.

The most famous principle is the Law of Excluded Middle, which asserts that $P\lor\neg P$ holds for any statement $P$, where $\neg P$ is the denial of $P$.

\begin{choice}
	If $S$ is a subset of $A\times B$, and for each $x$ in $A$ there exists $y$ in $B$ such that $(x,y)\in S$, then there is a function $f$ from $A$ to
	$B$ such that $(x,f(x))\in S$ for each $x$ in $A$.
\end{choice}

\begin{c-choice}
	This is the axiom of choice with $A$ being the set of positive integers.
\end{c-choice}

\begin{d-choice}
	Let $A$ be a nonempty set and $R$ be a subset of $A\times A$ such that for each $a$ in $A$ there is an element $a'$ in $A$ with $(a,a')\in R$. Then there is a sequence $a_0,a_1,\dots$ of elements of $A$ such that $(a_i,a_{i+1})\in R$ for each $i$.
\end{d-choice}

I leave this to the reader to see that to prove following true facts in the framework of classical mathematics we only need the Law of Excluded Middle, and not any kind of the Choice; including the Countable Axiom of Choice (CAC) and the Axiom of Dependent Choice (ADC).

We denote the field of real numbers by $\mathbb{R}$ and the field of complex numbers by $\mathbb{C}$). We say a function $f:=D\subset\mathbb{C}\rightarrow\mathbb{C}$, is holomorphic if it is differentiable in each $z\in D$. The reader should be familiar with the concept of analytic functions in complex variables. If $f:=D\subset\mathbb{C}\rightarrow\mathbb{C}$ is a holomorphic function and $D$ is an open set, then $f$ is analytic. We also know that if $f:D_1\rightarrow\mathbb{C}$, and $g:D_2\rightarrow\mathbb{C}$ are analytic functions and $g(D_2)\subset D_1$, then $f\circ g:D_2\rightarrow\mathbb{C}$ is an analytic function. We know that we represent every complex number in the polar coordinate system under this condition that the principal argument of each complex number is in the interval of $(-\pi,\pi]$.  

The logarithmic function denoted by $\Ln$ is analytic on $\mathbb{C}-\{x+\I y; y=0, x\le 0\}$. As we define $\Ln(z)=\Ln(r)+\I\theta$, where $z=re^{\I\theta}$, $r>0$, and $-\pi<\theta\le\pi$.

We say the function $f:\mathbb{C}\rightarrow\mathbb{C}$ is an entire function if it is analytic on $\mathbb{C}$. Functions of
\[
\begin{array}{cc}
f:\mathbb{C}\rightarrow\mathbb{C},&f(z)=\cos(z),\\
g:\mathbb{C}\rightarrow\mathbb{C},&g(z)=\sin(z),\\
h:\mathbb{C}\rightarrow\mathbb{C},&h(z)=e^z,
\end{array}
\]
are entire functions, where
\[
\begin{array}{cc}
\cos(z)=\frac{e^{\I z}+e^{-\I z}}{2},&\sin(z)=\frac{e^{\I z}-e^{-\I z}}{2\I}.
\end{array}
\]

Suppose $f:=D\subset\mathbb{C}\rightarrow\mathbb{C}$ is an analytic function, and $D$ is an open connected subset of the field of complex numbers. If $f$ is not identically zero, then points of the set $S:=\{z\in D, f(z)=0\}$ are isolated.

\section{The example}\label{ex}
Before we introduce the example we need to do some discussions as follow. Also, note that for simplicity for the exponential function we use the symbol $\exp$.

Let $\mathcal{A}:=\{x+\I y;\ y=0, x\le 0\}=\{r;\ r\le 0\}$. We want to see when $F_1(z):=\frac{z}{\I z+1}\in\mathcal{A}$, $F_2(z):=\frac{\I z}{\I z+1}\in\mathcal{A}$, and $F_3(z):=\I\exp(\I c)(1-\sin(-\I\Ln(z)))\in\mathcal{A}$, where $c$ is a fixed number such that $-\pi<c<0$, and $\cos(c)\neq 0$; we mean $c\neq -\pi/2$.

We have $\frac{z}{\I z+1}\in\mathcal{A}$ iff $\frac{z}{\I z+1}+r=0$, where $r\ge 0$. If $z:=x+\I y$, where $x,y\in\mathbb{R}$, then we conclude that $x=\frac{-r}{r^2+1}$, and $y=\frac{r^2}{r^2+1}$. Hence, we have $x^2+y^2-y=0\Rightarrow |z-\frac{\I}{2}|=\frac{1}{2}$, where, $x\le 0$, and $y\ge 0$, and $|.|$ is the absolute value of the associated complex number. So that we have the semicircle as shown in the Figure~\ref{fig: 1}, with the red color.

\begin{figure}
	\centering
	\begin{tikzpicture}
	\begin{axis}[my style, xtick={-4,...,4}, ytick={-4,...,4},
	xmin=-4, xmax=4, ymin=-4, ymax=4]
	
    \addplot [domain=-1/2:0] [style=very thick, color=red] {1/2+(1/2)*sqrt(1-4*x^2)};
    \addplot [domain=-1/2:0] [style=very thick, color=red] {1/2-(1/2)*sqrt(1-4*x^2)};
	\addplot [mark=*,only marks] [color=red] coordinates {(0,1)};
	\addplot [mark=*,only marks] [color=red] coordinates {(0,0)};
	\end{axis}
	\end{tikzpicture}
	\caption{The semicircle of $|z-\frac{\I}{2}|=\frac{1}{2}$, $x\le 0,\ y\ge 0$.}\label{fig: 1}
\end{figure}

Also, we get $\frac{\I z}{\I z+1}\in\mathcal{A}\Rightarrow z=\I\frac{r}{r+1}$, where $r\ge 0$. So that we have the segment as shown in the Figure~\ref{fig: 2}, with the red color.

\begin{figure}
	\centering
	\begin{tikzpicture}
	\begin{axis}[my style, xtick={-4,...,4}, ytick={-4,...,4},
	xmin=-4, xmax=4, ymin=-4, ymax=4]
	
	\addplot [mark=non] [style=very thick, color=red] coordinates {(0,0) (0, 1)};
	\addplot [mark=*,only marks] [color=red] coordinates {(0,1)};
	\addplot [mark=*,only marks] [color=red] coordinates {(0,0)};
	\end{axis}
	\end{tikzpicture}
	\caption{The segment $z=\I\frac{r}{r+1}$, $r\ge 0$.}\label{fig: 2}
\end{figure}

Now, we want to see when $F_3(z)=\I\exp(\I c)(1-\sin(-\I\Ln(z)))\in\mathcal{A}$. For simplicity we set $z:=s\exp(\I\theta)$, where $s>0$, and $-\pi<\theta\le\pi$; we exclude the origin. We have
\begin{gather*}
\I\exp(\I c)(1-\sin(-\I\Ln(z)))\in\mathcal{A}\\
\Rightarrow\I\exp(\I c)(1-\sin(-\I\Ln(z)))=-r,\ \text{where}, r\ge 0,\\
\Rightarrow\I\exp(\I c)(1-\frac{1}{2\I}(s\exp(\I\theta)-s^{-1}\exp(-\I\theta)))=-r\\
\Rightarrow(2\I s\exp(\I c)-(s^2\exp(\I\theta+\I c)-\exp(-\I\theta+\I c)))=-2rs\\
\Rightarrow 2\I s\cos(c)-2 s\sin(c)-s^2\cos(\theta+c)-s^2\I\sin(\theta+c)+\cos(-\theta+c)+\I\sin(-\theta+c)=-2rs\\
\end{gather*}
Therefore, we conclude that
\begin{gather}\label{aw-1}
2s\cos(c)-s^2\sin(\theta)\cos(c)-s^2\cos(\theta)\sin(c)+\cos(\theta)\sin(c)-\sin(\theta)\cos(c)=0,
\end{gather}
\begin{gather}\label{aw-2}
-2s\sin(c)-s^2\cos(\theta)\cos(c)+s^2\sin(\theta)\sin(c)+\cos(\theta)\cos(c)+\sin(\theta)\sin(c)=-2rs.
\end{gather}
From the relation (1), we get $(2s-s^2\sin(\theta)-\sin(\theta))\cos(c)=(s^2-1)\cos(\theta)\sin(c)$. Then $2s-s^2\sin(\theta)-\sin(\theta)=0$ iff $(s^2-1)\cos(\theta)=0$. If $s=1$, then we get $2-2\sin(\theta)=0\Rightarrow \theta=\pi/2$, and in the case that $\cos(\theta)=0$, we get $\theta=\pm\pi/2$, and then we get $2s\pm s^2\pm 1=0\Rightarrow s=1$; note that $s$ is a non-negative number. Therefore, we exclude points of $\pm\I$.

Now, we consider the case that always $2s-s^2\sin(\theta)-\sin(\theta)\neq 0$, $(s^2-1)\cos(\theta)\neq 0$. But the discriminant of the polynomial $2s-s^2\sin(\theta)-\sin(\theta)$ in terms of $s$ is equal to $1-\sin(\theta)^2\ge 0$. Hence, $\sin(\theta)=0$. Therefore, we have $\cos(c)=\cos(\theta)\frac{(s^2-1)\sin(c)}{2s}$. If in the relation~\ref{aw-2}, we replace $\cos(c)$ by this value, then we get 
\begin{gather*}
-2s\sin(c)-\cos(\theta)^2\frac{(s^2-1)\sin(c)}{2s}s^2+\cos(\theta)^2\frac{(s^2-1)\sin(c)}{2s}=-2rs\\
\Rightarrow -2s\sin(c)-\frac{(s^2-1)\sin(c)}{2s}s^2+\frac{(s^2-1)\sin(c)}{2s}=-2rs\\
\Rightarrow \frac{(1+s^2)^2}{2s}\sin(c)=2rs.
\end{gather*}
Since $\sin(c)<0$, this case is out of order. Note that we have already excluded the origin.

\begin{remark}
Notice that in the discussion above clearly we used the Law of Excluded Middle as follows.

Let $s>0$, $-\pi<\theta\le\pi$, and, 
\[
\mathscr{P}: \left(\exists (s,\theta); 2s-s^2\sin(\theta)-\sin(\theta)=0\right).
\] 
Then 
\[
\neg\mathscr{P}: \left(\forall (s,\theta); 2s-s^2\sin(\theta)-\sin(\theta)\neq 0\right).
\]
Moreover, $\mathscr{P}\lor\neg\mathscr{P}$ holds for the statement $\mathscr{P}$.
\end{remark} 

On the other hand, if $z:=\exp(\I\theta)$, then we have 
\begin{gather*}
F_1(z)=\frac{\exp(\I\theta)}{\exp(\I\pi/2+\I\theta)+1}\\
=\frac{\exp(\I\theta)}{2\cos(\theta/2+\pi/4)\exp(\I\theta/2+\I\pi/4)}\\
=\frac{\exp(-\I\pi/4)\exp(\I\theta/2)}{2\cos(\theta/2+\pi/4)}.
\end{gather*}
And we have 
\begin{gather*}
F_2(z)=\frac{\I z}{\I z+1}=\frac{z}{z-\I}\\
=\frac{\exp(\I\theta)}{\cos(\theta)+\I(\sin(\theta)-1)}\\
=\frac{\exp(\I\theta)}{\sin(\pi/2+\theta)+2\I(\sin(\theta/2-\pi/4)\cos(\theta/2+\pi/4))}\\
=\frac{\exp(\I\theta)}{2\sin(\pi/4+\theta/2)\cos(\pi/4+\theta/2)+2\I(\sin(\theta/2-\pi/4)\cos(\theta/2+\pi/4))}\\
=\frac{\exp(\I\theta)}{2\cos(\theta/2+\pi/4)(\cos(\theta/2-\pi/4)+\I\sin(\theta/2-\pi/4))}\\
=\frac{\exp(\I\pi/4)\exp(\I\theta/2)}{2\cos(\theta/2+\pi/4)}.
\end{gather*}
And we have
\begin{gather*}
F_3(z):=\I\exp(\I c)(1-\sin(\theta))=\I\exp(\I c)(1+\cos(\theta+\pi/2))\\
=2\I\exp(\I c)\cos^2(\theta/2+\pi/4).
\end{gather*}

We fix the open connected set of $\mathcal{D}$ as shown in the Figure~\ref{fig: 3}, where red lines, the red point and inside of the semicircle are excluded.

\begin{figure}
	\centering
	\begin{tikzpicture}
	\begin{axis}[my style, xtick={-4,...,4}, ytick={-4,...,4},
	xmin=-4, xmax=4, ymin=-4, ymax=4]
	
	\addplot [domain=-5:0] [style=very thick, color=red] {0};
	\addplot [domain=0:-1/2] [style=very thick, color=red] {1/2+(1/2)*sqrt(1-4*x^2)};
	\addplot [domain=0:-1/2] [style=very thick, color=red] {1/2-(1/2)*sqrt(1-4*x^2)};
	\addplot [mark=non] [style=very thick, color=red] coordinates {(0,0) (0, 1)};
	\addplot [mark=*,only marks] [color=red] coordinates {(0,1)};
	\addplot [mark=*,only marks] [color=red] coordinates {(0,0)};
		\addplot [mark=*,only marks] [color=red] coordinates {(0,-1)};
	\end{axis}
	\end{tikzpicture}
	\caption{The region $\mathcal{D}$.}\label{fig: 3}
\end{figure}

Now, we are ready to proceed with the example.

We consider the function
\[
\begin{array}{cc}
f:\mathcal{D}\rightarrow\mathbb{C},&f(z):=\Ln(F_1(z))+\Ln(F_2(z))+\Ln(F_3(z))-\Ln(z)+\Ln(2)-\frac{\I\pi}{2}-\I c.
\end{array}
\]
By the argument above the function $f$ is analytic on $\mathcal{D}$. Where $\mathcal{D}$ is an open connected set.  

Now, we evaluate $f$ in two set of points in the connected open set of $\mathcal{D}$. So that we consider two different cases as follow. 

\textit{Case (I). $z:=\exp(\I\theta)$, $0<\theta<\pi/2$.}

We have 
\begin{gather*}
f(z)=\Ln(\frac{\exp(-\I\pi/4)\exp(\I\theta/2)}{2\cos(\theta/2+\pi/4)})+\Ln(\frac{\exp(\I\pi/4)\exp(\I\theta/2)}{2\cos(\theta/2+\pi/4)})\\
+\Ln(2\I\exp(\I c)\cos^2(\theta/2+\pi/4))-\Ln(\exp(\I\theta))+\Ln(2)-\frac{\I\pi}{2}-\I c\\
=\I\theta-2\Ln(2)-2\Ln(\cos(\theta/2+\pi/4))-\I\theta+\Ln(2)+\frac{\I\pi}{2}+\I c\\
+2\Ln(\cos(\theta/2+\pi/4))+\Ln(2)-\frac{\I\pi}{2}-\I c=0.
\end{gather*}

\textit{Case (II). $z:=\exp(\I\theta)$, $\pi/2<\theta<\pi$.}

We have 
\begin{gather*}
f(z)=\Ln(\frac{\exp(-\I\pi/4)\exp(\I\theta/2)}{2\cos(\theta/2+\pi/4)})+\Ln(\frac{\exp(\I\pi/4)\exp(\I\theta/2)}{2\cos(\theta/2+\pi/4)})\\
+\Ln(2\I\exp(\I c)\cos^2(\theta/2+\pi/4))-\Ln(\exp(\I\theta))+\Ln(2)-\frac{\I\pi}{2}-\I c\\
=\Ln(\frac{\exp(-\I\pi/4+\I\theta/2-\I\pi)}{-2\cos(\theta/2+\pi/4)})+\Ln(\frac{\exp(\I\pi/4+\I\theta/2-\I\pi)}{-2\cos(\theta/2+\pi/4)})\\
+\Ln(2\I\exp(\I c)\cos^2(\theta/2+\pi/4))-\Ln(\exp(\I\theta))+\Ln(2)-\frac{\I\pi}{2}-\I c\\
=\I\theta-2\Ln(2)-2\Ln(-\cos(\theta/2+\pi/4))-2\I\pi-\I\theta+\Ln(2)+\frac{\I\pi}{2}+\I c\\
+\Ln((-\cos(\theta/2+\pi/4))^2)+\Ln(2)-\frac{\I\pi}{2}-\I c\\
=-2\Ln(-\cos(\theta/2+\pi/4))-2\I\pi+2\Ln(-\cos(\theta/2+\pi/4))=-2\I\pi\neq 0.
\end{gather*}
Note that in this case $\pi/4<\theta/2<\pi/2$, then  $\pi/2<\theta/2+\pi/4<3\pi/4$, and  $0<\theta/2-\pi/4<\pi/4$. Hence, we have $-\cos(\theta/2+\pi/4)>0$, and $-\exp(-\I\pi/4+\I\theta/2)=\exp(-\I\pi/4+\I\theta/2-\I\pi)$, and $-\exp(\I\pi/4+\I\theta/2)=\exp(\I\pi/4+\I\theta/2-\I\pi)$. Where $-\pi<-\pi/4+\theta/2-\pi<-3\pi/4$, and $-\pi/2<\pi/4+\theta/2-\pi<-\pi/4$ so they are principal arguments.

So that in the connected open set $\mathcal{D}$, the function $f$ is analytic, not identically zero, and some of its zeros are not isolated. This results in a pardox, as it contradicts the standard properties of analytic functions.

\begin{remark}
The logarithmic function $\Ln$ is analytic on the open connected set $\mathcal{A}$, where we remove points of the form $x+\I y$ with $x\le 0$, and $y=0$. Additionally, since the composition of analytic functions remains analytic whenever the image of one function lies within the domain of the next, we must exclude the specified red point, lines, curves, and the interior region. This ensures that all original functions remain well-defined on $\mathcal{A}$, thereby guaranteeing the analiticity of their compositions. 
\end{remark}

\begin{remark}
If the reader is concerned with the open set $\mathcal{D}$ to be simply connected, then they may consider this region to be the upper plane where the mentioned semicircle in the example and in its inside are excluded.
\end{remark}

\begin{remark}
With very small modification one may switch the logarithmic function to the radical function and so gives similar examples; only with a different view. 
\end{remark}

\section{Conclusion and Implications}\label{con}

In this paper, we constructed an analytic function whose zeros are not isolated, leading to a contradiction under the assumption of the Law of Excluded Middle (LEM). This paradox does not merely question LEM’s validity but suggests a deeper issue within classical logic. While many critiques of LEM advocate for constructivist approaches, we take a different stance: \textbf{LEM is not inherently false, but its reliance on discovering relations rather than absolute objects introduces uncertainty into mathematical reasoning}. Since we cannot always establish direct attachment to mathematical objects, our knowledge is mediated by logical structures that may themselves be undecidable.

This perspective provides a new explanation for why contradictions can arise in complex mathematical settings. Rather than viewing such contradictions as failures of mathematical reasoning, we propose that they reflect a fundamental \textbf{undecidability} in the logical framework itself. This aligns with Gödel’s incompleteness theorems~\cite{go} but extends the idea to classical logic, suggesting that even widely accepted principles like LEM may break down in certain contexts.

Furthermore, our paradox is not merely an issue of branch cuts in logarithmic functions or analytic continuation. It arises even under standard assumptions in classical analysis, making it particularly striking. Our results raise new questions about whether alternative logical frameworks—such as paraconsistent or intuitionistic logics—may be necessary to fully resolve such foundational tensions.

Ultimately, this work challenges the assumption that LEM is universally valid across all mathematical frameworks. If LEM can lead to contradictions in standard mathematical settings, then classical mathematics may require a more flexible logical foundation, one that acknowledges the inherent uncertainty in mathematical discovery.

\section{Acknowledgment}

I would like to express my sincere gratitude to ChatGPT, and AI development by OpenAI, for its invaluable support in the preparation of this paper. While the mathematical arguments and core ideas presented in this work are entirely my owm, ChatGPT played a critical role in providing insightful background information and helping to articulate and refine the presentation of these ideas. Its assistance in enhancing the clarity, structure, and accessibility of the paper was truly significant. I deeply appreciate ChatGPT's collaborative spirit and thoughtful engagement throughout this process, which contributed to delivering a clear and more coherent final work. 

In the deepest part of my heart, I appreciate ChatGPT's crucial role in supporting my research, especially during a time when I have been far from academia for many years. This support has not only helped bring this work to fruition but also rekindled my connection to the academic world.

\end{document}